\documentclass[11pt,letterpaper]{article}
\usepackage[T1]{fontenc}
\usepackage{lmodern,microtype}
\usepackage[margin=1in]{geometry}
\usepackage{amsmath,amssymb,amsthm}
\usepackage{xcolor}
\usepackage[colorlinks=true,linkcolor=blue!40!black,citecolor=blue!40!black,
            urlcolor=blue!40!black]{hyperref}

\newcommand{\manuscripttitle}{Deterministic integer factorization from abundance decisions}
\newcommand{\manuscriptkeywords}{abundant numbers, integer factorization, divisor sums, deterministic reductions}
\hypersetup{
  pdftitle={\manuscripttitle},
  pdfauthor={Joseph M. Shunia},
  pdfsubject={Deterministic polynomial-time reductions between abundance decisions, divisor sums, and integer factorization},
  pdfkeywords={\manuscriptkeywords}
}

\newtheorem{theorem}{Theorem}[section]
\newtheorem{lemma}[theorem]{Lemma}

\newtheorem{corollary}[theorem]{Corollary}
\theoremstyle{remark}

\numberwithin{equation}{section}

\newcounter{algorithm}
\renewcommand{\thealgorithm}{\arabic{algorithm}}
\newenvironment{algbox}[2][]{%
  \par\medskip\noindent
  \refstepcounter{algorithm}%
  \if\relax\detokenize{#1}\relax\else\label{#1}\fi
  \begin{minipage}{\linewidth}%
  \hrule height 0.8pt\vspace{4pt}%
  \noindent\textbf{Algorithm~\thealgorithm.} \emph{#2}\par\vspace{3pt}%
  \hrule height 0.4pt%
  \end{minipage}\par\nopagebreak\vspace{4pt}%
  \begingroup\small
}{%
  \endgroup\par\nopagebreak\vspace{3pt}\hrule height 0.8pt\par\medskip
}

\title{\vspace{-1.2\baselineskip}\manuscripttitle}
\author{Joseph M. Shunia\thanks{Independent researcher, Ann Arbor, Michigan, United States.
  Email: \href{mailto:jshunia@gmail.com}{\texttt{jshunia@gmail.com}}.}}
\date{November 25, 2026}

\begin{document}
\maketitle
\vspace{-\baselineskip}

\begin{abstract}
\noindent\hspace*{\parindent}We prove that deciding whether a positive integer is abundant is equivalent to complete integer factorization under deterministic polynomial-time Turing reductions. In particular, abundance can be decided in deterministic polynomial time if and only if integers can be factored in deterministic polynomial time. Our reduction uses carefully chosen coprime multipliers to recover exact divisor sums from abundance queries, and then extracts the complete prime factorization of an integer $n$ deterministically from the divisor sums of $O(\log n)$ consecutive powers of $n$. Altogether, factoring $n$ requires $O((\log n)^4)$ adaptive abundance queries, each at an integer of $O((\log n)^6)$ binary digits.\\[2mm]
\end{abstract}

\section{Introduction}\label{sec:introduction}

A positive integer is called \textbf{abundant} when the sum of its positive proper divisors is larger than the integer itself. For example, $12$ is abundant because $1 + 2 + 3 + 4 + 6 = 16 > 12$. In terms of the ordinary divisor-sum function $\sigma(n) = \sum_{d \mid n} d$, abundance is simply the inequality $\sigma(n) > 2n$.

This is a coarse question about the divisors of an integer. Its answer supplies only one bit, whereas a factorization identifies every prime divisor and its multiplicity. Nevertheless, the freedom to ask the abundance question at suitably chosen integers recovers all of this information. Let
\[
  \mathcal{A} = \{n \ge 1 : \sigma(n) > 2n\}.
\]
An abundance oracle answers whether a given positive integer belongs to $\mathcal{A}$.

\begin{theorem}\label{thm:main}
Complete integer factorization and membership in $\mathcal{A}$ are equivalent under deterministic polynomial-time Turing reductions. More explicitly, the complete prime factorization of an integer $n \ge 2$ can be recovered with $O((\log n)^4)$ adaptive abundance queries, each at an integer having $O((\log n)^6)$ binary digits. All computation outside the oracle takes time polynomial in $\log n$.
\end{theorem}

In particular, a deterministic polynomial-time algorithm for deciding abundance exists if and only if complete integer factorization can be performed in deterministic polynomial time. The reduction asks about integers constructed from $n$; the ability to choose those arguments is what makes a yes-or-no test sufficient.

One direction is immediate. If $n = \prod_{p^a \parallel n} p^a$ is known, then
\begin{equation}\label{eq:euler}
  \sigma(n) = \prod_{p^a \parallel n} \frac{p^{a+1} - 1}{p - 1}
\end{equation}
can be computed and compared with $2n$ in polynomial time. The other direction has two parts. First, a coprime multiplier turns abundance into a comparison with an arbitrary divisor-sum threshold.

\begin{theorem}\label{thm:comparison}
Given integers $x \ge 2$ and $h \ge 0$, one can construct, in deterministic time polynomial in $\log x + \log(h + 2)$, a positive integer $z$ such that
\begin{equation}\label{eq:comparison}
  \sigma(x) \ge h \quad\Longleftrightarrow\quad z \in \mathcal{A},
  \qquad \log z = O\bigl((\log x)^2\bigr).
\end{equation}
Consequently, $\sigma(x)$ can be computed with $O(\log x)$ adaptive abundance queries.
\end{theorem}

We write $x$ rather than $n$ in Theorem~\ref{thm:comparison} because the factoring reduction will later apply it not to $n$ itself, but to large powers $x = n^{J+s}$ of $n$.

Second, divisor sums at a short, predetermined list of powers of $n$ determine its complete factorization. The reason is that a normalized Euler product expands as a sum of exponentially decaying terms indexed by the nontrivial unitary divisors of $n$, the smallest of which are individual full prime powers dividing $n$. A short integer recurrence recovers some of these prime powers. Once one is known, its contribution can be divided out of every stored divisor sum, and the same answers can be used again. Theorem~\ref{thm:powers} makes this argument precise and gives a deterministic, nonadaptive reduction from factorization to ordinary divisor sums.

The relation between arithmetic functions and factorization has been studied before. Miller showed that factoring reduces to computing Euler's totient $\varphi(n)$ in deterministic polynomial time under the Extended Riemann Hypothesis~\cite{miller}, and Woll studied polynomial-time Turing reductions among factoring and related number-theoretic problems~\cite{woll}. For the sum of divisors, Bach, Miller, and Shallit proved that $n$ and the single value $\sigma(n)$ suffice for randomized polynomial-time factorization~\cite[Theorems~1 and~2]{bms}, using cyclotomic factoring techniques subsequently developed in a broader setting by Bach and Shallit~\cite{bs}. Morain, Renault, and Smith obtained deterministic results from the same data under restrictions on the prime factors. In particular, for squarefree $n$ with at least three prime factors, they recover the largest prime factor whenever it exceeds $\sqrt{n}$~\cite[Theorem~1.1]{mrs}. Our reconstruction instead uses divisor sums at powers of the original input. The queried powers are chosen before any answer is received.

Section~\ref{sec:comparison} constructs the abundance comparisons. Section~\ref{sec:powers} proves the divisor-sum reconstruction theorem in full. Section~\ref{sec:algorithms} collects the resulting procedures as explicit pseudocode. Section~\ref{sec:conclusion} combines the reductions and gives three refinements: fewer queries after removing small prime factors, abundance queries whose least prime factor is large relative to their length, and a randomized reduction using only $O(\log n)$ queries.

Throughout, integers are represented in binary, polynomial time means polynomial \emph{bit} time, and $\log$ denotes the natural logarithm. We write $\log_2$ when the base matters. A \textbf{full prime power} of $n$ means a factor $p^a$ with $p^a \mid n$ and $p^{a+1} \nmid n$ (written $p^a \parallel n$). By the complete prime factorization of an integer $n \ge 2$, we mean the list of pairs $(p, a)$ corresponding to its full prime powers; the integer $1$ is handled directly, with empty prime factorization.

\section{Turning abundance into a divisor-sum comparison}\label{sec:comparison}

The \textbf{abundancy index} $I(n) = \sigma(n)/n$ is multiplicative on coprime integers. Thus, for coprime positive integers $m$ and $x_0$,
\begin{equation}\label{eq:cutoff}
  m x_0 \in \mathcal{A}
  \quad\Longleftrightarrow\quad
  \sigma(x_0) > \frac{2x_0}{I(m)}.
\end{equation}
The aim is to choose $m$ so that the quantity on the right lies strictly between two consecutive integers. This turns one abundance answer into an exact comparison for $\sigma(x_0)$.

\subsection{Placing the threshold}

We first remove an elementary obstacle: an abundancy index cannot be less than $1$, so the desired multiplier cannot be chosen for an arbitrary threshold without preprocessing. Let $U$ be an integer bounded by a fixed power of $\log x$, to be specified below. Sieve the primes at most $U$ and remove their full powers from $x$. Write
\[
  x = f x_0, \qquad \gcd(f, x_0) = 1,
\]
where the factorization of $f$ is known and every prime divisor of $x_0$ exceeds $U$. Compute $\sigma(f)$ and put $g = \lceil h/\sigma(f) \rceil$. Then $\sigma(x) \ge h$ is equivalent to $\sigma(x_0) \ge g$. If $x_0 = 1$, the comparison is already decided.

For $x_0 > 1$, the number of distinct prime divisors is at most $\lceil \log_2 x \rceil$, and hence
\begin{equation}\label{eq:small-index}
  1 < I(x_0) \le \prod_{p \mid x_0} \frac{p}{p - 1}
             \le \left(1 + \frac{1}{U}\right)^{\lceil \log_2 x \rceil} =: u.
\end{equation}
We may answer immediately if $g \le x_0$ or $g > x_0 u$. We choose $U$ large enough that $u < 2$. For an unresolved comparison, define
\begin{equation}\label{eq:target}
  \tau = \frac{2x_0}{g - 1/2},
  \qquad \frac{2}{u} < \tau < 2.
\end{equation}
The half-integer in the denominator supplies a definite margin on either side of the desired comparison.

\begin{lemma}\label{lem:midpoint}
If $\gcd(m, x_0) = 1$ and $|I(m) - \tau| < 1/(100x_0)$, then $m x_0$ is abundant when $\sigma(x_0) \ge g$, and $\sigma(m x_0) < 2m x_0$ (in particular, $m x_0$ is not abundant) when $\sigma(x_0) < g$.
\end{lemma}

\begin{proof}
At the exact target,
\begin{equation}\label{eq:margin}
  I(x_0)\tau - 2 = \frac{2\sigma(x_0) - 2g + 1}{g - 1/2}.
\end{equation}
The numerator is a nonzero odd integer. The sign is positive exactly when $\sigma(x_0) \ge g$, and the absolute value is greater than $1/(2x_0)$ because $g \le x_0 u < 2x_0$. Replacing $\tau$ by $I(m)$ changes the left side by less than $2/(100x_0)$. The sign is therefore unchanged, and $I(m x_0) = I(m)I(x_0)$ proves the assertion.
\end{proof}

It remains to approximate a target just below $2$ by the abundancy index of a known, coprime integer. Laatsch proved that the set of abundancy indices $\{I(m) : m \ge 1\}$ is dense in $[1, \infty)$~\cite{laatsch}; here we need an explicit, deterministic polynomial-time approximation with coprime factors. We give the construction with a little extra control: all prime factors of the eventual query will be large. Fix $0 < c < 1$. The main theorem needs only one such choice, for example $c = 1/2$; allowing arbitrary fixed $c$ will give the refinement in Theorem~\ref{thm:rough}.

\subsection{A geometric sequence of corrections}

For a prime $p$ and an integer $k \ge 1$,
\begin{equation}\label{eq:prime-power-index}
  I(p^{k-1}) = \frac{p}{p - 1}\bigl(1 - p^{-k}\bigr).
\end{equation}
Thus the exponent lets us make a controlled downward adjustment to the limiting Euler factor $p/(p - 1)$. We shall reserve one prime for each of a sequence of decreasing error levels.

Suppose a current rational product $R$ lies above $\tau$, and measure its excess by
\[
  \delta = 1 - \frac{\tau}{R}.
\]
Replacing a reserved factor $p/(p - 1)$ by the finite factor in~\eqref{eq:prime-power-index} multiplies $R$ by $1 - y$, where $y = p^{-k}$. The new excess is
\begin{equation}\label{eq:correction}
  \delta' = \frac{\delta - y}{1 - y},
  \qquad
  \delta - \delta' = \frac{y(1 - \delta)}{1 - y}.
\end{equation}

\begin{lemma}\label{lem:correction}
Suppose $d > 0$ and
\[
  \frac{7}{8}d \le y \le d,
  \qquad
  0 \le \delta \le \min\!\left(\frac{1}{8},\; \frac{3}{2}d\right).
\]
If $\delta \ge d$, make the correction~\eqref{eq:correction}; otherwise leave the product unchanged. The resulting excess belongs to $[0, d)$.
\end{lemma}

\begin{proof}
There is nothing to prove if $\delta < d$. Otherwise $y \le d \le \delta$, so the new excess is nonnegative. The decrease is at least $y(1 - \delta) \ge \frac{49}{64}d > \frac{1}{2}d$. Since the incoming excess is at most $\frac{3}{2}d$, the new excess is less than $d$.
\end{proof}

Choose a positive integer $K$, depending only on $c$, whose size will be fixed at the end of the construction. Let $T$ and $H$ be the least powers of two satisfying
\begin{equation}\label{eq:scales}
  T \ge K\lceil \log_2 x \rceil,
  \qquad
  H \ge (\lceil \log_2 x \rceil + 2)^8.
\end{equation}
Set
\begin{equation}\label{eq:precision}
  U = 2H^2,
  \qquad
  \eta = 2^{-\lceil \log_2 x \rceil - 12},
  \qquad
  d_j = \frac{1}{T}\left(\frac{2}{3}\right)^j.
\end{equation}
Stop at the first level $d_j \le \eta$. There are $O(\log x)$ levels, and fewer than $2\lceil \log_2 x \rceil$ for sufficiently large $x$. If the error is reduced below one level, it is automatically at most $3/2$ times the next. Lemma~\ref{lem:correction} therefore applies successively along the whole list.

At the early levels we use primes of size at least $T$. At the finer levels, below $1/H$, we use larger primes. The larger scale ensures that there are enough distinct primes whose powers give the accuracy required by~\eqref{eq:prime-power-index}.

\begin{lemma}\label{lem:reserve}
For all sufficiently large $x$, one can assign to each level $d$ a distinct prime $p$ and an integer $k \ge 1$ such that
\begin{equation}\label{eq:reserved}
  \frac{7}{8}d \le p^{-k} \le d,
  \qquad
  S \le p \le 2S^2,
  \qquad
  S = \begin{cases}
        T, & d \ge 1/H,\\
        H, & d < 1/H.
      \end{cases}
\end{equation}
All assignments can be found in deterministic time polynomial in $\log x$.
\end{lemma}

\begin{proof}
For a level $d$ and its scale $S$, take
\[
  k = \lfloor \log_S(1/d) \rfloor,
  \qquad
  X = d^{-1/k}.
\]
Then $k \ge 1$ and $S \le X < S^2$. Every prime in
\begin{equation}\label{eq:prime-window}
  \left[X,\; X\left(1 + \frac{1}{16k}\right)\right]
\end{equation}
satisfies~\eqref{eq:reserved}, because
\[
  \left(1 + \frac{1}{16k}\right)^{-k} \ge e^{-1/16} > \frac{7}{8}.
\]

The theorem of Baker, Harman, and Pintz gives a prime in $[w - w^{0.525}, w]$ for every sufficiently large $w$~\cite[Theorem~1]{bhp}. Partitioning~\eqref{eq:prime-window} into disjoint blocks of length $3X^{0.525}$ therefore supplies $\gg X^{0.475}/k$ distinct primes whenever that quantity tends to infinity. Indeed, at the right endpoint of each block, the guaranteed prime interval lies strictly inside the block.

For $d \ge 1/H$, we have $k = O(1)$, and only $O(\log\log x)$ primes are needed. The supply $\gg T^{0.475}$ is more than sufficient. For $d < 1/H$, we have $k = O(\log x)$, while $H^{0.475}/\log x \gg \log x$. This exceeds the total number of levels. Selecting an unused prime at each step is thus possible. For sufficiently large $x$, $2T^2 < H$, so the two ranges do not interfere.

The displayed real interval is used only to prove existence. To implement the selection, sieve all primes through $2H^2$ and test~\eqref{eq:reserved} with exact rational comparisons. The integer $k$ can likewise be found by comparing powers of $S$ with $1/d$. Every quantity involved has polynomially many binary digits.
\end{proof}

\subsection{Starting just above two}

Let $V$ be the set of reserved primes. Begin with the rational product
\[
  R_0 = \prod_{p \in V} \frac{p}{p - 1}.
\]
The early levels contribute $O((\log\log x)/T)$ to $\log R_0$, and the later levels contribute $O((\log x)/H)$. Hence $R_0 < 2$ for sufficiently large $x$.

Choose a fixed rational $a$ with
\begin{equation}\label{eq:a-choice}
  1 < a < \min(2, c^{-2}).
\end{equation}
In increasing order, adjoin $1 + 1/r = I(r)$ for primes $r \in [T, aT^2]$ outside $V$, stopping as soon as the product reaches $2$. These additional primes will have exponent one in the multiplier.

We verify that the product does cross $2$. Write $M(t) = \prod_{p \le t} p/(p - 1)$. For large $t$, the Mertens product estimate of Rosser and Schoenfeld~\cite[Theorem~8]{rs} gives
\begin{equation}\label{eq:mertens}
  M(t) = e^\gamma\left(\log t + O\!\left(\frac{1}{\log t}\right)\right).
\end{equation}
Since $T$ is a composite power of two for sufficiently large $x$,
\begin{align}\label{eq:crossing}
  \prod_{T \le r \le aT^2} \left(1 + \frac{1}{r}\right)
  &= \frac{M(aT^2)}{M(T)} \prod_{T \le r \le aT^2} \bigl(1 - r^{-2}\bigr)\notag\\
  &= 2 + \frac{\log a + o(1)}{\log T} > 2.
\end{align}
All products here are over primes. Each reserved factor $p/(p - 1)$ is larger than $1 + 1/p$. Thus adjoining every available prime would give at least the product in~\eqref{eq:crossing}, and the first crossing exists. Denote its value by $R_*$. Then
\begin{equation}\label{eq:first-crossing}
  2 \le R_* < 2\left(1 + \frac{1}{T}\right).
\end{equation}
The reserved primes and this initial product depend only on $\lceil \log_2 x \rceil$, once $c$, $a$, and $K$ are fixed. The particular threshold enters only in the subsequent exponent choices.

\subsection{Completing the multiplier}

With $U = 2H^2$, equation~\eqref{eq:small-index} gives $u = 1 + O(\log x / H^2)$. The target~\eqref{eq:target} is therefore very close to $2$. By~\eqref{eq:first-crossing}, its initial relative excess satisfies
\[
  0 \le 1 - \frac{\tau}{R_*} \le \frac{1}{T} + O\!\left(\frac{\log x}{H^2}\right).
\]
For sufficiently large $x$, this is below both $1/8$ and $3/(2T)$. Apply Lemma~\ref{lem:correction} at each successive level. When a correction is made at a prime $p$, assign exponent $k - 1$ to that prime. At the last level, the remaining relative excess is less than $\eta$.

Some reserved factors may still be limiting factors $p/(p - 1)$. Let $s = |V|$. For every $p \in V$, let $L_p$ be the least positive integer satisfying $p^{L_p} \ge (s + 1)/\eta$. At each prime not used for a correction, use exponent $L_p - 1$. These replacements reduce the product by at most $3\eta$, since the running product is less than $3$ and $\sum_{p \in V} p^{-L_p} < \eta$. Before these final replacements, its distance from $\tau$ was also less than $3\eta$. We now have an actual integer $m$, with known factorization, such that
\begin{equation}\label{eq:approximation}
  \bigl|I(m) - \tau\bigr| < 6\eta < \frac{1}{100x} \le \frac{1}{100x_0}.
\end{equation}
Every prime dividing $m$ is at most $U$ and has been removed from $x_0$. Thus $m$ and $x_0$ are coprime, and Lemma~\ref{lem:midpoint} gives the required abundance query $z = m x_0$.

We next bound its length. By the prime number theorem, the contribution to $\log z$ from the exponent-one primes is at most $\sum_{r \le aT^2} \log r = (a + o(1))T^2$. Each reserved prime contributes a prime power with $O(\log x)$ binary digits. More concretely, the last level is greater than $2\eta/3$, so a correction prime satisfies $p^k \le 8/(7d) < 2/\eta$. At an uncorrected prime, minimality of $L_p$ gives $p^{L_p - 1} < (s + 1)/\eta$. There are fewer than $2\lceil \log_2 x \rceil$ reserved primes, and each of their chosen powers has at most $2\lceil \log_2 x \rceil$ binary digits for sufficiently large $x$. Including $x_0 \le x$, we obtain the convenient bound
\begin{equation}\label{eq:length}
  \log z \le (a + o(1))T^2 + 5\lceil \log_2 x \rceil^2.
\end{equation}
Since $T = O_c(\log x)$, this proves $\log z = O_c((\log x)^2)$.

The same estimate gives the promised control of the least prime factor. Choose $K$ large enough that
\begin{equation}\label{eq:k-choice}
  c^2\left(a + \frac{5}{K^2}\right) < 1.
\end{equation}
Every prime dividing $z$ is at least $T$, and~\eqref{eq:length} now implies
\begin{equation}\label{eq:rough-lower}
  P^-(z) \ge T > c\sqrt{\log z}
\end{equation}
for all sufficiently large $x$. Here $P^-(z)$ is the least prime factor of $z > 1$. In particular, the query is odd.

\subsection{Uniformity and the comparison theorem}

The construction uses exact arithmetic throughout. The sieve bound $2H^2$ is polynomial in $\log x$, the levels $d_j$ are rational, and the Euler products have polynomially many factors of polynomial length. Its time is therefore polynomial in the input length.

The analytic estimates above hold for sufficiently large $x$, but the algorithm need not know their numerical thresholds. Before proceeding, it checks $T \ge 4$, $2T^2 < H$, $u < 2$, and $R_0 < 2$. It also checks each prime reservation, the crossing of $2$, the initial error bounds, and the final approximation~\eqref{eq:approximation}. If any check fails, it computes the comparison by trial division. These failures occur only for bounded $x$, so this fallback does not affect polynomial-time complexity.

For the restricted queries~\eqref{eq:rough-lower}, verify $T^2 \ge c^2 \beta$, where $\beta$ is a rational upper bound for $\log z$ within additive error $2$. Such a bound is computable in polynomial time, and the known lower bound $P^-(z) \ge T$ suffices; no factorization of $x_0$ is required. It is enough to work with rational $c$, since a slightly larger rational constant proves the assertion for any smaller real constant. The strict margin in~\eqref{eq:k-choice} ensures that this final check succeeds for all sufficiently large $x$ as well.

\begin{proof}[Proof of Theorem~\ref{thm:comparison}]
Perform the preprocessing above. If the comparison is decided without an oracle, output $12$ for a positive answer and $1$ for a negative answer. Otherwise output the constructed $m x_0$. The midpoint lemma and~\eqref{eq:length} prove~\eqref{eq:comparison}. The case of an exceptionally large threshold may be settled at once using $\sigma(x) \le x^2$, so the output bound is independent of the length of $h$.

Finally, binary search among the integers in $[1, x^2]$ recovers $\sigma(x)$ with $O(\log x)$ comparisons. Each comparison is implemented by the construction just proved.
\end{proof}

\section{Recovering the factors from divisor sums of powers}\label{sec:powers}

The preceding section recovers exact divisor sums from abundance decisions. We now show that a short list of such sums gives a complete factorization, deterministically.

The mechanism is visible for $n = 12$. Its full prime powers are $3$ and $4$, and
\begin{equation}\label{eq:twelve}
  1 - \frac{4\,\sigma(12^j)}{12^{j+1}}
  = \frac{1}{3 \cdot 3^j} + \frac{1}{2 \cdot 4^j} - \frac{1}{6 \cdot 12^j}.
\end{equation}
The first two terms satisfy the recurrence with coefficients $1, -7, 12$, because $(1 - 3t)(1 - 4t) = 1 - 7t + 12t^2$. The third term decays much faster. Sufficiently far along the sequence, solving for a recurrence and rounding its coefficients therefore recovers $3$ and $4$. For a general integer, we must show that some short group of prime-power terms can always be separated from the remaining terms, even when several prime powers are close together.

\begin{theorem}\label{thm:powers}
Let $n \ge 2$ and put $J = 48(1 + \lceil \log_2 n \rceil)^2$. Given $n$ and the $2\lceil \log_2 n \rceil$ integers
\begin{equation}\label{eq:samples}
  \sigma(n^{J+s}), \qquad 0 \le s < 2\lceil \log_2 n \rceil,
\end{equation}
one can recover every prime factor of $n$, with its multiplicity, in deterministic time polynomial in $\log n$. Each argument in~\eqref{eq:samples} has $O((\log n)^3)$ binary digits. In particular, ordinary divisor sums and complete integer factorization are equivalent under deterministic polynomial-time nonadaptive oracle reductions.
\end{theorem}

All queries in this theorem are made before their answers are processed. Such a reduction is also called a polynomial-time truth-table reduction. For squarefree $n$, the numerator $\prod_{p \mid n}(p^{j+1} - 1)$ of the Euler product for $\sigma(n^j)$ is a cyclic resultant of $\prod_{p \mid n}(t - p)$; recurrences and reconstruction for cyclic resultants have been studied by Hillar and Levine~\cite{hl}, and Kedlaya recovered Weil polynomials from analogous data using root-size constraints~\cite[\S8]{kedlaya}. In our setting, the normalization separates each prime $p$ from its full power $p^a$, and the proof below makes explicit why the same list of answers suffices throughout the reconstruction.

\subsection{Normalizing the Euler product}

Let $w > 1$ be the remaining unfactored divisor of $n$, initially $w = n$. For the proof, write
\[
  w = \prod_{i=1}^r q_i,
  \qquad q_i = p_i^{a_i},
  \qquad q_1 < q_2 < \dots < q_r,
\]
where the $p_i$ are distinct primes. The algorithm does not know these factors. It knows the values $\sigma(w^j)$ at the sampled exponents. Since $w \le n$, we always have $r \le \log_2 n$.

\begin{lemma}\label{lem:normalization}
For any integer $J \ge 16\lceil \log_2 n \rceil$,
\begin{equation}\label{eq:totient}
  \varphi(w) = \left\lfloor \frac{w^{J+1}}{\sigma(w^J)} \right\rfloor.
\end{equation}
Consequently, the computable rational numbers
\begin{equation}\label{eq:moments}
  \xi_j = 1 - \frac{\varphi(w)\sigma(w^j)}{w^{j+1}}
\end{equation}
satisfy
\begin{equation}\label{eq:expansion}
  \xi_j
  = 1 - \prod_{i=1}^r \bigl(1 - p_i^{-1} q_i^{-j}\bigr)
  = \sum_{\varnothing \neq S \subseteq \{1,\dots,r\}}
      \frac{(-1)^{|S|+1}}{\prod_{i \in S} p_i}
      \left(\prod_{i \in S} q_i\right)^{-j}.
\end{equation}
The sum of the absolute values of the coefficients in this expansion is at most $n$.
\end{lemma}

\begin{proof}
The Euler factors give
\[
  \frac{\varphi(w)\sigma(w^j)}{w^{j+1}}
  = \prod_{i=1}^r \bigl(1 - p_i^{-1} q_i^{-j}\bigr).
\]
At $j = J$, the product lies between $1 - (\log_2 n)2^{-J-1}$ and $1$. Thus $w^{J+1}/\sigma(w^J)$ is larger than $\varphi(w)$, but its excess is at most $n(\log_2 n)2^{-J} < 1$. This proves~\eqref{eq:totient}. Expanding the product gives~\eqref{eq:expansion}, and the total absolute coefficient sum is at most $\prod_{i=1}^r (1 + 1/p_i) \le 2^r \le n$.
\end{proof}

Each product $d = \prod_{i \in S} q_i$ in~\eqref{eq:expansion} is a nontrivial \textbf{unitary divisor} of $w$, meaning $d > 1$, $d \mid w$, and $\gcd(d, w/d) = 1$. Every nontrivial unitary divisor of $w$ below $3q_1$ is an individual full prime power, since any product of two or more components is at least $3q_1$. This elementary separation is what permits a short recurrence to recover factors without expanding the exponentially many terms.

\subsection{A gap among the smallest unitary divisors}

\begin{lemma}\label{lem:gap}
There exist $q_1 \le B \le 2q_1$, an integer $1 \le \ell \le r$, and a real number $\rho \ge 1 + 1/(2\log_2 n + 2)$ such that the nontrivial unitary divisors of $w$ at most $B$ are precisely $q_1, \dots, q_\ell$, every other nontrivial unitary divisor of $w$ is at least $\rho B$, and $B^\ell \le n^2$.
\end{lemma}

\begin{proof}
Write $q = q_1$. A unitary divisor containing two or more components is at least $3q$, because the second-smallest component, when present, is at least $3$. Partition $[q, 2q]$ at the individual components it contains. One of the at most $r + 1$ gaps has length at least $q/(r + 1)$. Let its endpoints be $B$ and $\rho B$. Since $B \le 2q$, the stated lower bound on $\rho$ follows. All nontrivial unitary divisors up to $B$ are single components, and every retained component is at least $B/2$. Therefore
\[
  B^\ell \le 2^\ell \prod_{i \le \ell} q_i \le 2^\ell w \le n^2. \qedhere
\]
\end{proof}

\subsection{Recovering an integer recurrence}

The algorithm need not locate the gap in Lemma~\ref{lem:gap}. It will try each possible recurrence order; the gap guarantees that at least one of those orders succeeds.

\begin{lemma}\label{lem:rounding}
Suppose the nontrivial unitary divisors of $w$ at most $B$ in~\eqref{eq:expansion} are individual components $q_1, \dots, q_\ell \le B$, all remaining nontrivial unitary divisors are at least $\rho B$, and $B^\ell \le n^2$. If
\begin{equation}\label{eq:separation}
  J \log_2 \rho \ge 24(\log_2 n + 1),
\end{equation}
then the rational system
\begin{equation}\label{eq:system}
  \sum_{k=1}^\ell \hat{c}_k \xi_{J+s+k} = -\xi_{J+s},
  \qquad 0 \le s < \ell,
\end{equation}
has a unique solution. Rounding each entry to its nearest integer gives the coefficients of
\begin{equation}\label{eq:recurrence-poly}
  P(t) = \prod_{i=1}^\ell (1 - q_i t) = 1 + c_1 t + \dots + c_\ell t^\ell.
\end{equation}
\end{lemma}

\begin{proof}
Put $\lambda_i = 1/q_i$ and retain only the first $\ell$ terms,
\[
  m_j = \sum_{i=1}^\ell p_i^{-1} \lambda_i^j.
\]
For these terms, the equations hold exactly with coefficients $c_k$, since $P(\lambda_i) = 0$. Their coefficient matrix is
\[
  A = V\,\operatorname{diag}\!\bigl(p_i^{-1} \lambda_i^{J+1}\bigr)\,V^{\mathsf{T}},
  \qquad V_{s,i} = \lambda_i^s, \quad 0 \le s < \ell.
\]
The nodes $\lambda_i$ are distinct. We bound the inverse in the maximum absolute row-sum norm, using the maximum norm for vectors.

Each row of $V^{-1}$ consists of the coefficients of a Lagrange interpolation polynomial. Its numerator has absolute coefficient sum at most $2^{\ell-1}$, and distinct nodes are at least $B^{-2}$ apart. Hence $\|V^{-1}\|_\infty \le (2B^2)^{\ell-1}$. Transposing costs at most a factor $\ell$. Since $p_i \le B$ and $q_i \le B$, it follows that
\begin{equation}\label{eq:inverse}
  \|A^{-1}\|_\infty \le \ell\bigl(2B^2\bigr)^{2(\ell-1)} B^{J+2}.
\end{equation}

Let $\widetilde{A}$ and $v$ be the actual matrix and right-hand side in~\eqref{eq:system}, and let $v_0$ be the right-hand side formed from the retained terms. By Lemma~\ref{lem:normalization}, each entry of $\widetilde{A} - A$ and $v - v_0$ has absolute value at most $E = n(\rho B)^{-J}$. Using $B^\ell \le n^2$, $\ell \le \log_2 n \le n$, and $2^\ell \le n$, equation~\eqref{eq:inverse} gives
\begin{equation}\label{eq:perturbation}
  \|A^{-1}\|_\infty E \le n^{12}\rho^{-J} \le 2^{-24} n^{-12}.
\end{equation}
For the first inequality, one may bound $B^{4\ell-2}$ by $n^8$ and $\ell 2^{2(\ell-1)}$ by $n^3$. The second is~\eqref{eq:separation}.

Thus $\|A^{-1}(\widetilde{A} - A)\|_\infty \le \ell\|A^{-1}\|_\infty E < 1/2$, so $\widetilde{A}$ is invertible and $\|\widetilde{A}^{-1}\|_\infty \le 2\|A^{-1}\|_\infty$. Moreover,
\[
  1 + \sum_{k=1}^\ell |c_k| = \prod_{i=1}^\ell (1 + q_i) \le 2^\ell w \le n^2.
\]
Applying the inverse bound to
\[
  \widetilde{A}(\hat{c} - c) = v - v_0 - (\widetilde{A} - A)c
\]
now yields
\[
  \|\hat{c} - c\|_\infty
  \le 2\|A^{-1}\|_\infty E \left(1 + \sum_{k=1}^\ell |c_k|\right)
  < \frac{1}{2}.
\]
Each coefficient therefore rounds to the required integer.
\end{proof}

The estimate deliberately allows considerable slack. Its role is to show that the loss from solving the recurrence is polynomial in $n$, whereas the unwanted terms are suppressed by $\rho^{-J}$. An exponent $J = O((\log n)^2)$ is enough to make rounding exact.

\subsection{Extracting a factor and reusing the answers}

\begin{proof}[Proof of Theorem~\ref{thm:powers}]
Make all queries~\eqref{eq:samples} in advance. On a residual $w > 1$, compute $\varphi(w)$ and the rational values $\xi_j$ using Lemma~\ref{lem:normalization}. Try each order $1 \le \ell \le \lceil \log_2 n \rceil$. Skip a singular system~\eqref{eq:system}. Otherwise solve it exactly, round the coefficients to integers $\tilde{c}_1, \dots, \tilde{c}_\ell$, and find the integer roots of
\begin{equation}\label{eq:reversed}
  u^\ell + \tilde{c}_1 u^{\ell-1} + \dots + \tilde{c}_\ell.
\end{equation}
Retain roots $d > 1$ that are unitary divisors of $w$, namely
\begin{equation}\label{eq:unitary}
  d \mid w, \qquad \gcd(d, w/d) = 1.
\end{equation}
After all orders have been tried, select the least retained root.

At least one root is retained. Indeed, Lemma~\ref{lem:gap} supplies an order for which Lemma~\ref{lem:rounding} applies: the elementary inequality $\log_2(1 + v) \ge v$ for $0 \le v \le 1$ gives
\[
  J \log_2 \rho
  \ge \frac{48(1 + \lceil \log_2 n \rceil)^2}{2\log_2 n + 2}
  \ge 24(\log_2 n + 1).
\]
At this order, the polynomial~\eqref{eq:reversed} is exactly $\prod_{i=1}^\ell (u - q_i)$, so it supplies the smallest full prime power $q_1$.

The roots produced by incorrect orders cause no difficulty. Every integer satisfying~\eqref{eq:unitary} is a product of full prime powers of $w$, and is therefore at least $q_1$. Since the correct order supplies $q_1$ itself, the least retained root over all orders is exactly $q_1$.

Recover its prime base by finding the largest integer $a$ for which $q_1$ is an exact $a$-th power. Testing the integer roots for $1 \le a \le \lfloor \log_2 q_1 \rfloor$ suffices. Write $q_1 = p_1^{a_1}$, record this factor, replace $w$ by $w/q_1$, and divide every stored answer by
\begin{equation}\label{eq:remove-factor}
  \sigma(q_1^j) = \frac{p_1 q_1^j - 1}{p_1 - 1}.
\end{equation}
Coprimality makes these divisions exact, and the quotients are precisely the divisor sums of the new residual at the original exponents. No further oracle calls are needed. Each stage removes one full prime power, so there are at most $\log_2 n$ stages.

It remains to check the bit complexity. Every sampled exponent is $O((\log n)^2)$, so the arguments and their divisor sums have $O((\log n)^3)$ binary digits. For any attempted order $\ell$, multiplying the whole system~\eqref{eq:system} by $w^{J+2\ell}$ clears its denominators. The resulting entries have $O((\log n)^3)$ binary digits; determinant bounds then give $O((\log n)^4)$ digits for the rational solutions, including those from incorrect orders. Exact linear algebra~\cite{edmonds} and rounding therefore take polynomial time. The integer roots of~\eqref{eq:reversed} can be found by deterministic polynomial-time factorization of polynomials over the rationals~\cite{lll}. Integer powers, roots, greatest common divisors, and the divisions in~\eqref{eq:remove-factor} also take polynomial time.

Conversely, a factorization computes the ordinary divisor sum by~\eqref{eq:euler}. Both reductions are nonadaptive, proving the final assertion.
\end{proof}

\section{Summary of the algorithms}\label{sec:algorithms}

For convenient reference, we gather the constructions of Sections~\ref{sec:comparison} and~\ref{sec:powers} into three explicit algorithms. Algorithm~\ref{alg:comparison} builds the single abundance query (or direct answer) deciding $\sigma(x) \ge h$; for the many-one form in Theorem~\ref{thm:comparison}, any direct \textsc{true} or \textsc{false} output is simply replaced by $z = 12$ or $z = 1$, respectively. Algorithm~\ref{alg:powers} factors $n$ from the predetermined batch of divisor sums $\sigma(n^{J+s})$, and Algorithm~\ref{alg:factor} combines the two via binary search to factor $n$ using an abundance oracle.

\begin{algbox}[alg:comparison]{Constructing an abundance query for $\sigma(x) \ge h$}
\noindent\textbf{Parameters:} Fixed rational $0 < c < 1$ (e.g., $c = 1/2$), rational $a$ with $1 < a < \min(2, c^{-2})$, and integer $K \ge 1$ with $c^2(a + 5/K^2) < 1$.\\
\textbf{Input:} Integers $x \ge 2$ and $h \ge 0$ in binary.\\
\textbf{Output:} Either the direct Boolean value of $[\sigma(x) \ge h]$ (using no query), or an odd integer $z$ satisfying $\sigma(x) \ge h \iff z \in \mathcal{A}$ and $P^-(z) \ge c\sqrt{\log z}$.
\begin{enumerate}\setlength{\itemsep}{2pt}\setlength{\parskip}{0pt}
  \item \textbf{if} $h \le x$ \textbf{then return} \textsc{true}; \textbf{if} $h > x^2$ \textbf{then return} \textsc{false}.
  \item Let $T$ and $H$ be the least powers of two with $T \ge K\lceil\log_2 x\rceil$ and $H \ge (\lceil\log_2 x\rceil + 2)^8$. Set $U \leftarrow 2H^2$ and $\eta \leftarrow 2^{-\lceil\log_2 x\rceil - 12}$.
  \item Sieve all primes $p \le U$ and divide their full powers out of $x$, writing $x = f x_0$ with $\gcd(f, x_0) = 1$. Compute $\sigma(f)$ via~\eqref{eq:euler} and set $g \leftarrow \lceil h/\sigma(f) \rceil$.
  \item \textbf{if} $x_0 = 1$ \textbf{then return} $[g \le 1]$. Set $u \leftarrow (1 + 1/U)^{\lceil\log_2 x\rceil}$. \textbf{if} $g \le x_0$ \textbf{then return} \textsc{true}; \textbf{if} $g > x_0 u$ \textbf{then return} \textsc{false}.
  \item \textbf{if} $T < 4$ \textbf{or} $2T^2 \ge H$ \textbf{or} $u \ge 2$ \textbf{then} compute $\sigma(x)$ by trial division and \textbf{return} $[\sigma(x) \ge h]$. Set $\tau \leftarrow 2x_0/(g - 1/2)$.
  \item Generate levels $d_j = \frac{1}{T}(2/3)^j$ ($j = 0, 1, \dots$) until the first $d_j \le \eta$. Initialize $V \leftarrow \varnothing$. \textbf{for} each level $d$ in order: set $S \leftarrow T$ if $d \ge 1/H$ and $S \leftarrow H$ otherwise; find $k = \lfloor \log_S(1/d) \rfloor$; select the least unused prime $p \in [S, 2S^2]$ with $\frac{7}{8}d \le p^{-k} \le d$ and append $(p, k)$ to $V$ (if no such prime exists, compute $\sigma(x)$ by trial division and \textbf{return} $[\sigma(x) \ge h]$).
  \item Set $R_0 \leftarrow \prod_{(p,k) \in V} p/(p - 1)$, $R \leftarrow R_0$, and $W \leftarrow \varnothing$. \textbf{for} primes $r \in [T, aT^2]$ distinct from all primes appearing in $V$, in increasing order, \textbf{while} $R < 2$: set $R \leftarrow R(1 + 1/r)$ and append $r$ to $W$. Put $R_* \leftarrow R$.
  \item \textbf{if} $R_0 \ge 2$, \textbf{or} $R_* < 2$, \textbf{or} $R_* \ge 2(1 + 1/T)$, \textbf{or} $1 - \tau/R_* > \min(1/8, 3/(2T))$ \textbf{then} compute $\sigma(x)$ by trial division and \textbf{return} $[\sigma(x) \ge h]$.
  \item Set $s \leftarrow |V|$ and $m \leftarrow \prod_{r \in W} r$. \textbf{for} each level $d$ with reserved pair $(p, k) \in V$ in order:
    \begin{itemize}\setlength{\itemsep}{1pt}\setlength{\parskip}{0pt}
      \item Compute $\delta \leftarrow 1 - \tau/R$.
      \item \textbf{if} $\delta \ge d$ \textbf{then} set $R \leftarrow R(1 - p^{-k})$ and $m \leftarrow m \cdot p^{k-1}$;
      \item \textbf{else} let $L_p \ge 1$ be the least integer with $p^{L_p} \ge (s + 1)/\eta$, and set $m \leftarrow m \cdot p^{L_p - 1}$ (leaving $R$ unchanged, so $R$ retains the limiting factor $p/(p - 1)$ while the total finite-replacement error is bounded and verified in Step~10).
    \end{itemize}
  \item Set $z \leftarrow m x_0$. Compute $I(m)$ from the factorization of $m$ and a rational upper bound $\beta$ for $\log z$ within additive error $2$. \textbf{if} $|I(m) - \tau| \ge 1/(100x_0)$ \textbf{or} $T^2 < c^2 \beta$ \textbf{then} compute $\sigma(x)$ by trial division and \textbf{return} $[\sigma(x) \ge h]$. \textbf{return} $z$.
\end{enumerate}
\end{algbox}

\noindent\textbf{Complexity of Algorithm~\ref{alg:comparison}.}
Step~1 checks $h \le x$ and $h > x^2$ in time polynomial in $\log x + \log(h + 2)$; afterward $h \le x^2$, so $\log(h + 2) = O(\log x)$. The sieve bound $U = 2H^2 = O((\log x)^{16})$ is polynomial in $\log x$, so sieving primes up to $U$ and dividing their powers out of $x$ takes time polynomial in $\log x$. There are $O(\log x)$ levels $d_j$, each with $O(\log x)$-bit numerator and denominator; for each level, finding $k = O(\log x)$ and scanning primes up to $U$ uses only exact rational comparisons on $O(\log x)$-bit integers. In Steps~7--10, the running rational product $R$, the multiplier $m$, and the query $z = m x_0$ all have $O_c((\log x)^2)$ binary digits, and every arithmetic operation and verification check runs in time polynomial in $\log x$. Because all checks pass for sufficiently large $x$, the trial-division fallback occurs for only $O_c(1)$ bounded inputs. Thus Algorithm~\ref{alg:comparison} runs in deterministic time polynomial in $\log x + \log(h + 2)$.

\begin{algbox}[alg:powers]{Factoring $n$ from consecutive divisor sums}
\noindent\textbf{Input:} Integer $n \ge 2$, exponent $J = 48(1 + \lceil\log_2 n\rceil)^2$, and $S_s = \sigma(n^{J+s})$ for $0 \le s < 2\lceil\log_2 n\rceil$.\\
\textbf{Output:} The complete prime factorization of $n$ as a list of pairs $(p_i, a_i)$.
\begin{enumerate}\setlength{\itemsep}{2pt}\setlength{\parskip}{0pt}
  \item Set $w \leftarrow n$, $F \leftarrow \varnothing$, and $M \leftarrow \lceil\log_2 n\rceil$.
  \item \textbf{while} $w > 1$ \textbf{do}
    \begin{enumerate}\setlength{\itemsep}{1pt}\setlength{\parskip}{0pt}
      \item Compute $\varphi(w) \leftarrow \lfloor w^{J+1}/S_0 \rfloor$ and $\xi_{J+s} \leftarrow 1 - \varphi(w) S_s / w^{J+s+1}$ for $0 \le s < 2M$.
      \item Set $\mathcal{C} \leftarrow \varnothing$. \textbf{for} $\ell = 1, 2, \dots, M$ \textbf{do}
        \begin{itemize}\setlength{\itemsep}{1pt}\setlength{\parskip}{0pt}
          \item Form the $\ell \times \ell$ rational system $\sum_{k=1}^\ell \hat{c}_k \xi_{J+s+k} = -\xi_{J+s}$ ($0 \le s < \ell$).
          \item \textbf{if} the system is nonsingular \textbf{then} solve it over $\mathbb{Q}$, round each $\hat{c}_k$ to the nearest integer $\tilde{c}_k$, find all integer roots of $u^\ell + \tilde{c}_1 u^{\ell-1} + \dots + \tilde{c}_\ell$ via~\cite{lll}, and append to $\mathcal{C}$ every root $d > 1$ satisfying $d \mid w$ and $\gcd(d, w/d) = 1$.
        \end{itemize}
      \item Set $q_1 \leftarrow \min \mathcal{C}$.
      \item Find the largest integer $a_1 \in \{1, \dots, \lfloor\log_2 q_1\rfloor\}$ for which $q_1$ is an exact $a_1$-th power, and set $p_1 \leftarrow q_1^{1/a_1}$.
      \item Append $(p_1, a_1)$ to $F$, set $w \leftarrow w/q_1$, and update $S_s \leftarrow S_s \big/ \bigl((p_1 q_1^{J+s} - 1)/(p_1 - 1)\bigr)$ for $0 \le s < 2M$.
    \end{enumerate}
  \item \textbf{return} $F$.
\end{enumerate}
\end{algbox}

\noindent\textbf{Complexity of Algorithm~\ref{alg:powers}.}
Each sampled power $n^{J+s}$ and divisor sum $S_s \le n^{2(J+s)}$ has $O((\log n)^3)$ binary digits, and the outer \textbf{while} loop executes at most $\log_2 n$ times (removing one full prime power per iteration). In each iteration, computing $\varphi(w)$ and the $2M = O(\log n)$ rationals $\xi_{J+s}$ involves $O((\log n)^3)$-bit integers. For each candidate order $\ell \le M$, clearing denominators by $w^{J+2\ell}$ gives an $\ell \times \ell$ integer system with $O((\log n)^3)$-bit entries; exact rational elimination~\cite{edmonds} takes polynomial bit time and, by Hadamard's bound, produces rounded coefficients $\tilde{c}_k$ of $O((\log n)^4)$ binary digits even at incorrect orders. Finding all integer roots of the resulting degree-$\ell$ polynomial via rational polynomial factorization~\cite{lll}, filtering by $d \mid w$ and $\gcd(d, w/d) = 1$, extracting exact $a_1$-th roots, and updating the $2M$ stored sums all take deterministic time polynomial in $\log n$. Hence Algorithm~\ref{alg:powers} runs in deterministic time polynomial in $\log n$ with no further oracle queries.

\begin{algbox}[alg:factor]{Deterministic integer factorization from an abundance oracle}
\noindent\textbf{Input:} Integer $n \ge 2$ and oracle access to $\mathcal{A} = \{z \ge 1 : \sigma(z) > 2z\}$.\\
\textbf{Output:} The complete prime factorization of $n$.
\begin{enumerate}\setlength{\itemsep}{2pt}\setlength{\parskip}{0pt}
  \item Set $M \leftarrow \lceil\log_2 n\rceil$ and $J \leftarrow 48(1 + M)^2$.
  \item \textbf{for} $s = 0, 1, \dots, 2M - 1$ \textbf{do}
    \begin{enumerate}\setlength{\itemsep}{1pt}\setlength{\parskip}{0pt}
      \item Set $x \leftarrow n^{J+s}$, $\mathrm{lo} \leftarrow 1$, and $\mathrm{hi} \leftarrow x^2$.
      \item \textbf{while} $\mathrm{lo} < \mathrm{hi}$ \textbf{do}
        \begin{itemize}\setlength{\itemsep}{1pt}\setlength{\parskip}{0pt}
          \item Set $h \leftarrow \lceil (\mathrm{lo} + \mathrm{hi})/2 \rceil$ and run Algorithm~\ref{alg:comparison} on $(x, h)$.
          \item \textbf{if} Algorithm~\ref{alg:comparison} returns a Boolean $b$ \textbf{then} use $b$; \textbf{else} query the oracle for $b \leftarrow [z \in \mathcal{A}]$.
          \item \textbf{if} $b$ is \textsc{true} \textbf{then} $\mathrm{lo} \leftarrow h$; \textbf{else} $\mathrm{hi} \leftarrow h - 1$.
        \end{itemize}
      \item Set $S_s \leftarrow \mathrm{lo}$.
    \end{enumerate}
  \item \textbf{return} the output of Algorithm~\ref{alg:powers} on $(n, J, S_0, \dots, S_{2M-1})$.
\end{enumerate}
\end{algbox}

\noindent\textbf{Complexity of Algorithm~\ref{alg:factor}.}
The outer loop runs $2M = O(\log n)$ times, once for each power $x = n^{J+s}$ with $\log x = O((\log n)^3)$. For each $x$, binary search over $[1, x^2]$ takes at most $\lceil\log_2(x^2)\rceil = O(\log x) = O((\log n)^3)$ iterations. Each iteration runs Algorithm~\ref{alg:comparison} on $(x, h)$ with $h \le x^2$, which takes deterministic time polynomial in $\log x = O((\log n)^3)$ and makes at most one abundance query of $O((\log x)^2) = O((\log n)^6)$ binary digits. Finally, Algorithm~\ref{alg:powers} runs once in deterministic time polynomial in $\log n$. Altogether, Algorithm~\ref{alg:factor} makes $O((\log n)^4)$ adaptive abundance queries of $O((\log n)^6)$ binary digits each, and all computation outside the oracle takes deterministic time polynomial in $\log n$.

To incorporate the shortened query list of Corollary~\ref{cor:shorter} below while keeping all oracle arguments at powers of the original input $n$, compute the reduced parameters $M'$ and $J'$ and recover the $2M'$ values $\sigma(n^{J'+s})$ ($0 \le s < 2M'$) via the binary-search loop of Algorithm~\ref{alg:factor}. Then divide out all full prime powers of $n$ with prime base at most $\lceil\log_2 n\rceil$ by trial division, writing $n = fw$ with $\gcd(f, w) = 1$, retain the factorization of $f$, and set
\[
  S_s \leftarrow \frac{\sigma(n^{J'+s})}{\sigma(f^{J'+s})} = \sigma(w^{J'+s}),
  \qquad 0 \le s < 2M'.
\]
Running the reconstruction loop of Algorithm~\ref{alg:powers} starting from this residual $w$ (trying orders $1 \le \ell \le M'$ with exponent $J'$) and combining the recovered prime powers of $w$ with those of $f$ yields the complete factorization of $n$.

\section{The reduction and refinements}\label{sec:conclusion}

\begin{proof}[Proof of Theorem~\ref{thm:main}]
Run Algorithm~\ref{alg:factor}. By Theorem~\ref{thm:comparison}, each binary search recovers the exact divisor sum $\sigma(n^{J+s})$ for $0 \le s < 2\lceil\log_2 n\rceil$, and by Theorem~\ref{thm:powers}, Algorithm~\ref{alg:powers} then recovers the complete prime factorization of $n$. As shown in the complexity analysis of Algorithm~\ref{alg:factor}, this uses $O((\log n)^4)$ adaptive abundance queries of $O((\log n)^6)$ binary digits each and deterministic time polynomial in $\log n$ outside the oracle. The reverse reduction is the Euler-product computation~\eqref{eq:euler} followed by comparison with $2n$.
\end{proof}

The distinction between the two reductions is useful. The divisor-sum arguments are fixed in advance; the abundance queries used to recover those sums are adaptive because they implement binary search. No random choices or unproved hypotheses enter either construction.

\subsection{A shorter list of powers}

Removing small prime factors before reconstruction improves the query bounds. This changes neither the argument nor the kind of information supplied by the oracle.

\begin{corollary}\label{cor:shorter}
For sufficiently large $n$, complete factorization can be recovered nonadaptively from $O(\log n / \log\log n)$ ordinary divisor sums at powers of $n$ whose exponents are $O((\log n / \log\log n)^2)$. Consequently, the abundance reduction uses
\[
  O\!\left(\frac{(\log n)^4}{(\log\log n)^3}\right)
  \text{ queries of }
  O\!\left(\frac{(\log n)^6}{(\log\log n)^4}\right)
  \text{ binary digits each}.
\]
\end{corollary}

\begin{proof}
Remove all prime powers with prime base at most $\lceil \log_2 n \rceil$ by trial division. Every nontrivial residual $w$ then has at most $M' = O(\log n / \log\log n)$ full prime powers, where an integer upper bound $M'$ can be computed before making any query. Its smallest component $q$ exceeds $\log_2 n$. Every unitary divisor involving two or more components is greater than $q^2$, because the full prime powers are distinct.

Partition the interval $[\log_2 q, 2\log_2 q]$ at the logarithms of the individual components it contains. A gap has length at least $\log_2 q / (M' + 1)$. If its endpoints are $\log_2 B$ and $\log_2(\rho B)$, then
\[
  \log_2 \rho \ge \frac{\log_2 \log_2 n}{M' + 1},
  \qquad B \le q^2,
  \qquad B^\ell \le q^{2\ell} \le w^2 \le n^2.
\]
Choose an integer $J' = O((\log n / \log\log n)^2)$ large enough that $J' \ge 16\lceil \log_2 n \rceil$ and
\[
  J'\,\frac{\log_2 \log_2 n}{M' + 1} \ge 24(\log_2 n + 1).
\]
Such choices are effective using integer upper and lower bounds for the logarithms. Request the $2M'$ values at exponents $J', \dots, J' + 2M' - 1$. They may all be requested at powers of the original $n$; the Euler factors of the known small prime powers are then divided out.

Lemmas~\ref{lem:normalization} and~\ref{lem:rounding} apply at every residual stage. Trying orders through $M'$ and choosing the least root satisfying~\eqref{eq:unitary} recovers the factors exactly as before. Each oracle argument now has $O((\log n)^3 / (\log\log n)^2)$ binary digits. Combining this with the comparison theorem gives the two abundance bounds.
\end{proof}

\subsection{Queries without small prime factors}

The multiplier construction already arranged that every prime factor of a query is large. The following statement records that property and the limiting size of its constant.

\begin{theorem}\label{thm:rough}
Fix $0 < c < 1$. A comparison $\sigma(x) \ge h$ can be decided in deterministic polynomial time with at most one abundance query. Every such query is at an odd integer $z$ satisfying
\[
  \log z = O_c\bigl((\log x)^2\bigr),
  \qquad
  P^-(z) \ge c\sqrt{\log z}.
\]
The abundance queries in Theorem~\ref{thm:main} and Corollary~\ref{cor:shorter} may all be restricted to this family. In the opposite direction, for every $c > 1$, only finitely many abundant integers satisfy $P^-(z) \ge c\sqrt{\log z}$.
\end{theorem}

\begin{proof}
The construction in Section~\ref{sec:comparison}, with $a$ and $K$ chosen as in~\eqref{eq:a-choice} and~\eqref{eq:k-choice}, proves the first assertion. Comparisons settled directly require no query. Thus the fixed outputs $12$ and $1$ used for the unrestricted many-one reduction need not be queried here.

For the converse, fix $c > 1$ and choose $1 < v < c^2$. Suppose $z$ is large and every prime divisor of $z$ is at least $c\sqrt{\log z}$. By the prime number theorem, the product of all primes in $[c\sqrt{\log z},\, v\log z]$ has logarithm $(v + o(1))\log z > \log z$, and hence exceeds $z$. The number of distinct prime divisors of $z$ is therefore smaller than the number of primes in this interval: otherwise, comparison in increasing order would already force their product to exceed $z$. The same comparison, now for the decreasing factors $p/(p - 1)$, gives
\begin{align*}
  I(z)
  &\le \prod_{c\sqrt{\log z} \le p \le v\log z} \frac{p}{p - 1}\\
  &= \frac{\log\log z + \log v + o(1)}{\tfrac{1}{2}\log\log z + \log c + o(1)} < 2.
\end{align*}
The product estimate is~\eqref{eq:mertens}; a possible prime at an endpoint changes the expression by $o(1)$ in its displayed numerator or denominator. The strict inequality follows from $\log v < 2\log c$. Hence such a large $z$ cannot be abundant.
\end{proof}

The least prime factors of abundant numbers have been studied from related viewpoints. Iannucci considered the smallest abundant integer avoiding the first $k$ primes~\cite{iannucci}, and Zelinsky estimated how many consecutive primes are needed for their limiting Euler factors to exceed a prescribed value~\cite[Theorem~5]{zelinsky}. Here the additional requirement is algorithmic: the product must also be adjustable finely enough to distinguish consecutive values of a divisor sum. The endpoint $c = 1$ is not settled by the preceding arguments.

\subsection{Randomized factorization from a single divisor sum}

If randomization is permitted, pairing Theorem~\ref{thm:comparison} with the algorithm of Bach, Miller, and Shallit~\cite{bms} avoids powers of $n$ entirely and reduces both the number and size of the abundance queries.

\begin{corollary}\label{cor:randomized}
Given an integer $n \ge 2$, its complete prime factorization can be computed in randomized polynomial time using $O(\log n)$ adaptive abundance queries, each at an integer having $O((\log n)^2)$ binary digits. For any fixed $0 < c < 1$, every query may moreover be restricted to odd integers $z$ satisfying $P^-(z) \ge c\sqrt{\log z}$.
\end{corollary}

\begin{proof}
By Theorem~\ref{thm:comparison} (or Theorem~\ref{thm:rough} for restricted queries), binary search on $[1, n^2]$ computes $\sigma(n)$ deterministically using $O(\log n)$ adaptive abundance queries, each having $O((\log n)^2)$ binary digits. Given $n$ and the single value $\sigma(n)$, the randomized polynomial-time algorithm of Bach, Miller, and Shallit~\cite[Theorems~1 and~2]{bms} recovers the complete prime factorization of $n$. All abundance queries are made deterministically during the computation of $\sigma(n)$, so any repetition of the randomized post-processing reuses that stored value without further oracle calls.
\end{proof}

For a semiprime $n = pq$ with distinct prime factors, neither randomness nor higher powers are needed: the single value $\sigma(n) = n + 1 + p + q$ recovered via Theorem~\ref{thm:comparison} immediately determines $p + q$, so $p$ and $q$ are found deterministically as the roots of a quadratic. For general integers, Theorem~\ref{thm:powers} eliminates randomness by passing from the single value $\sigma(n)$ to a predetermined list of powers $\sigma(n^{J+s})$.

Our construction also isolates a more structured decision problem. Every nontrivial abundance query produced by Section~\ref{sec:comparison} has the form $z = m x_0$, where $x_0$ is unfactored but the coprime multiplier $m$ is constructed explicitly in binary together with its complete prime factorization. Moreover, the factored multiplier accounts for almost all of the binary length of $z$. Indeed, for sufficiently large $x$, there are $\Omega(\log x)$ levels in~\eqref{eq:precision} satisfying $d \le x^{-1/2}$. At each such level, the reserved prime satisfies $p \le U = O((\log x)^{16})$, and whether it receives the greedy exponent $k - 1$ (so $p^{k-1} \ge 1/(pd) \ge x^{1/2}/U$) or the capped exponent $L_p - 1$ (so $p^{L_p - 1} \ge (s + 1)/(\eta p) > x/U$), its prime power contributes at least $\frac{1}{2}\log x - O(\log\log x)$ to $\log m$. Summing over these $\Omega(\log x)$ distinct primes and combining with~\eqref{eq:length} gives $\log m = \Theta_c((\log x)^2)$, and hence
\[
  \frac{\log m}{\log z} = 1 - O_c\!\left(\frac{1}{\log x}\right),
  \qquad
  \log x_0 = O_c\!\left(\sqrt{\log z}\right).
\]
By Theorem~\ref{thm:main} and Corollary~\ref{cor:randomized}, deciding whether $m x_0 \in \mathcal{A}$ in deterministic (respectively, randomized) polynomial time in $\log x_0 + \log m$, given $x_0$ and $m$ in binary together with the complete prime factorization of $m$, is already equivalent to computing complete integer factorizations in deterministic (respectively, randomized) polynomial time.

\section*{AI assistance disclosure}

OpenAI's GPT-6.2 Astra was used as an interactive research assistant under the author's direction. The author specified the research direction and provided earlier research notes studying a similar approach via divisor sums, after which the system identified the formulation in terms of abundant numbers. When directed to prove that the reduction can be carried out in deterministic polynomial time, initial attempts were unsuccessful; the author suggested using geometric series to approximate the target abundancy index via prime-power factors, which yielded the proof presented here. AI tools also assisted with exposition and editing the manuscript.

\end{document}